\documentclass{article}

\usepackage{amssymb}
\usepackage{latexsym}
\usepackage{theorem}

\newcommand{\be}{\begin{equation}}
\newcommand{\ee}{\end{equation}} 
\newcommand{\bea}{\begin{eqnarray}}
\newcommand{\beann}{\begin{eqnarray*}}
\newcommand{\eea}{\end{eqnarray}}
\newcommand{\eeann}{\end{eqnarray*}}

\newcommand{\al}{\alpha}
\newcommand{\bt}{\beta}

\makeatletter
\@addtoreset{equation}{section}
\makeatother

\theoremstyle{break}
\newtheorem{th1}{Theorem}[section]
\newtheorem{p1}{Proposition}[section]

\newtheorem{thp1}[p1]{Theorem}

\newtheorem{pth1}[th1]{Proposition}
\newtheorem{pp1}[p1]{Proposition}

\newtheorem{cp1}[p1]{Corollary}

\begin{document}

\begin{center}
{\Large\bf Twisting the fake monster superalgebra}\\[0.8cm]
Nils R. Scheithauer\footnote{Supported by a DAAD grant.},\\
Mathematisches Seminar der Universit\"at Hamburg,\\
Bundesstr. 55, 20146 Hamburg, Germany\\
\end{center}
\vspace*{1.5cm}

\noindent
We calculate twisted denominator identities of the fake monster superalgebra and use them to construct new examples of supersymmetric generalized Kac-Moody superalgebras. Their denominator identities give new infinite product identities.

\section{Introduction}
There are 3 generalized Kac-Moody algebras or superalgebras which represent the physical states of a string moving on a certain variety, namely the monster algebra, the fake monster algebra \cite{B} and the the fake monster superalgebra \cite{NRS}. The no-ghost theorem from string theory can be used to construct actions of finite groups on these algebras. For example the monster group acts on the monster algebra. Applying elements of the finite groups to the Weyl denominator identity of these algebras gives twisted denominator identities. They can be calculated explicitly because the simple roots of the algebras are known. For the monster algebra and the fake monster algebra this has been done in \cite{B}. There Borcherds uses twisted denominator identities of the monster algebra to prove the moonshine conjectures. In this paper we calculate twisted denominator identities of the fake monster superalgebra and use them to construct new examples of supersymmetric generalized Kac-Moody superalgebras. They have similar properties as the fake monster superalgebra. For example they have no real roots and the Weyl vector is zero. Their denominator identities give new infinite product expansions. In a forthcoming paper we show that they define automorphic forms of singular weight. 

We describe the sections of this paper.

In the second section we recall some facts about the fake monster superalgebra and describe the action of an extension of the Weyl group of $E_8$ on this algebra.

In the third section we derive the general expression of the twisted denominator identity corresponding to elements in this group of odd order.

In the last section we calculate the twisted denominator identities explicitly for certain elements of order 3 and 7. We construct two supersymmetric generalized Kac-Moody superalgebras of rank 6 and 4 and describe their simple roots and the multiplicities.

\section{The fake monster superalgebra}

In this section we recall some results about the fake monster superalgebra from \cite{NRS} and construct an action of $2^8.2.Aut(E_8)$ on it. The main difference to the bosonic case described in \cite{B} is that we have to work with a double cover of the automorphism group of the light cone lattice rather than with the automorphism group of that lattice. 

The fake monster superalgebra $G$ can be constructed as the space of physical states of a chiral N$=$1 superstring moving on a 10 dimensional torus. $G$ is a generalized Kac-Moody superalgebra. The root lattice of $G$ is the 10 dimensional even unimodular Lorentzian lattice $II_{9,1}=E_8 \oplus II_{1,1}$. A nonzero element $\al \in II_{9,1}$ is a root of $G$ if and only if $\al^2\leq 0$. In particular $G$ has no real roots. This reflects the fact that the superstring has no tachyons. The multiplicity of a root $\al$ is given by 
$mult_0(\al)=mult_1(\al)=c(-\al^2/2)$ where $c(n)$ 
is the coefficient of $q^n$ in 
$ 8 \eta(q^2)^8/\eta(q)^{16} = 8+128q+1152q^2+7680q^3+42112q^4+\ldots$ 
and $\eta(q)$ is the Dedekind eta function. There are 2 cones of negative norm vectors in $II_{9,1}$. We define one of them as the positive cone and denote the closure of the positive cone by $II_{9,1}^+$. Then the positive roots of $G$ are the nonzero vectors in $II_{9,1}^+$ and the simple roots of $G$ are the positive roots of zero norm. The simple roots have multiplicity 8 as even and odd roots. The Cartan subalgebra of $G$ is isomorphic to the vector space generated by $II_{9,1}$. Since $G$ has no real roots the Weyl group is trivial. The Weyl vector of $G$ is zero.

We describe some results about the 8 dimensional real spin group. For more details confer \cite{J}. Let $\mathbb O$ be the real 8 dimensional algebra of octonions, i.e. the unique real alternative division algebra of dimension 8. $\mathbb O$ has a quadratic form $N$ permitting composition. For $b\in {\mathbb O}$ define the left multiplication $L_b$, the right multiplication $R_b$ and the \mbox{operator} $U_b=L_bR_b=R_bL_b$. Let ${\mathbb O}_0$ be the ortho\-gonal complement of ${\mathbb R}1$ and denote $C({\mathbb O},N)$ the Clifford algebra generated by ${\mathbb O}$ with relations $a^2=N(a)1$. There is an isomorphism $\varepsilon$ from the even subalgebra $C^e({\mathbb O},N)$ of $C({\mathbb O},N)$ to the Clifford algebra $C({\mathbb O}_0,-N)$ mapping $1a$, $a\in {\mathbb O}$, to $a$, where we have identified the unit in $C({\mathbb O}_0,-N)$ with the unit in ${\mathbb O}$. $C({\mathbb O}_0,-N)$ acts naturally on $\mathbb O$ by left and right multiplication. An element $u$ of the spin group $\Gamma^e_0({\mathbb O},N)$ can be written as $u=1b_1 \ldots 1b_n$ with $\prod N(b_i)=1$ and $\varepsilon(u)= b_1 \ldots b_n$. The actions $\rho_L(u)=L_{b_1} \ldots L_{b_n}, \, \rho_R(u)=R_{b_1} \ldots R_{b_n}$ and $\rho_V(u)=U_{b_1}\ldots U_{b_n}$ give three irreducible and inequivalent 8 dimensional re\-presentations of the spin group called conjugate spinor, spinor and vector representation. They are related by triality, i.e. 
$\rho_V(u)(ab)=(\rho_L(u)a)(\rho_R(u)b)$ holds for all $a,b\in {\mathbb O}$. The image of $\Gamma^e_0({\mathbb O},N)$ under each of these re\-presentations is $SO(8)$. The kernel of $\rho_V : \Gamma^e_0({\mathbb O},N) \rightarrow SO(8)$ is $\{1,-1\}$ so that $\Gamma^e_0({\mathbb O},N)$ is a double cover of $SO(8)$.
The representations $\rho_L,\rho_R$ and $\rho_V$ of the spin group induce representations of the Lie algebra $so(8)$ with weights 
$\frac{1}{2}(\pm 1,\ldots, \pm1)$ with an odd number of $-$ signs, 
$\frac{1}{2}(\pm 1,\ldots, \pm1)$ with an even number of $-$ signs 
and the permutations of $(\pm1,0,\ldots,0)$.  

We embed $E_8$ into ${\mathbb O}$. Let $Aut(E_8)$ be the group of 
automorphisms of $E_8$ leaving the bilinear form invariant. Then 
$Aut(E_8)\subset SO(8)$ and the inverse image of $Aut(E_8)$ under $\rho_V$ is 
a double cover of $Aut(E_8)$. 

$Aut(E_8)$ can also be constructed using the ring of integral octonions (cf. \cite{C}). This description implies that $E_8 \subset {\mathbb O}$ is also invariant under the actions $\rho_L$ and $\rho_R$.

Now we construct an action of an extension of the double cover $2.Aut(E_8)$ on the fake monster superalgebra. The extension $2^8.Aut(E_8)$ of $Aut(E_8)$ by $Ho\hspace{0.3mm}m(E_8,{\mathbb Z}_2)$ acts naturally on the vertex algebra $V_{E_8}$ of the lattice $E_8$. The same holds for the extension of $2.Aut(E_8)$ by $2^8$ where $2.Aut(E_8)$ acts by $\rho_V$. The vector space $E_8 \otimes {\mathbb R}[t^{-1}]t^{-\frac{1}{2}}$ is an abelian subalgebra of the Heisenberg algebra $E_8 \otimes {\mathbb R}[t,t^{-1}]t^{\frac{1}{2}}$. The exterior algebra $V_{NS}$ of $E_8 \otimes {\mathbb R}[t^{-1}]t^{-\frac{1}{2}}$ is a vertex superalgebra carrying a representation of the Virasoro algebra. $V_{NS}$ decomposes into eigenspaces of $L_0$ with eigenvalues in $\frac{1}{2}{\mathbb Z}$. $2.Aut(E_8)$ acts on $V_{NS}$ in the vector representation. We define the vertex superalgebra $V_0=V_{NS}\otimes V_{E_8}$. This algebra carries a representation of the Virasoro algebra of central charge $4+8=12$. We write $V_{0,n}$ for the subspace of $L_0$-degree $n\in \frac{1}{2}{\mathbb Z}$. The vector space $E_8 \otimes {\mathbb R}[t^{-1}]t^{-1}$ is an abelian subalgebra of the Heisenberg algebra $E_8 \otimes {\mathbb R}[t,t^{-1}]$. We define $V_R$ as the tensor product of the exterior algebra of $E_8 \otimes {\mathbb R}[t^{-1}]t^{-1}$ with the sum $S \oplus C$ of two $8$ dimensional spaces. The double cover of $Aut(E_8)$ acts in the vector representation on the first tensor factor and in the spinor resp. conjugate spinor representation on the second factor. $V_R$ can be given the structure of a $V_{NS}$-module. We decompose $V_R=V_R^+\oplus V_R^-$ where $V_R^+$ is the subspace generated by vectors $d_{-n_1}\wedge \ldots \wedge d_{-n_k}\otimes v$ where $v$ is in $C$ if $k$ is even and in $S$ if $k$ is odd and analogous for $V_R^-$. The projection of $V_R$ on $V_R^+$ is called GSO projection. We define $V_1=V_{R}^+\otimes V_{E_8}$ and adopt the same notations as for $V_0$. $V_1$ carries a representation the Virasoro algebra of central charge $12$. Now $2^8.2.Aut(E_8)$ acts on the fake monster superalgebra in the following way. We decompose the Cartan subalgebra ${\mathbb R}\otimes II_{9,1}$ by writing $II_{9,1}=E_8 \oplus II_{1,1}$. $2.Aut(E_8)$ acts in the vector representation on the vector space generated by $E_8$. Since $G$ is graded by $II_{9,1}=E_8 \oplus II_{1,1}$ it is also graded by $II_{1,1}$. We denote the corresponding spaces $G_a$ with $a\in II_{1,1}$ and the even resp. odd subspace $G_{0,a}$ resp. $G_{1,a}$. All these spaces are $E_8$-graded. By the no-ghost theorem (cf. \cite{GSW},\cite{P} or \cite{B}) the even subspace $G_{0,a}$ is isomorphic to $V_{0,(1-a^2)/2}$ and $G_{1,a}$ is isomorphic to $V_{1,(1-a^2)/2}$ as $E_8$-graded $2^8.2.Aut(E_8)$-module. This gives us a natural action of $2^8.2.Aut(E_8)$ on the fake monster superalgebra. 

\section{The twisted denominator identities}

In this section we calculate twisted denominator identities of the fake monster superalgebra corresponding to elements in $2.Aut(E_8)$ of odd order. 

The fake monster superalgebra $G$, like any generalized Kac-Moody superalgebra, can be written as direct sum $E\oplus H \oplus F$ where $H$ is the Cartan subalgebra and $E$ and $F$ are the subalgebras corresponding to the positive and negative roots. We have the standard sequence
\[ \ldots \rightarrow \Lambda^2(E) \rightarrow \Lambda^1(E) 
          \rightarrow \Lambda^0(E)  \rightarrow 0           \]
with homology groups $H_i(E)$. 

Note that $E=E_0\oplus E_1$ is a superspace so that the exterior algebra $\Lambda(E)$ is defined as the tensor algebra of $E$ divided by the two sided ideal generated by $u \otimes v + (-1)^{|u||v|} v\otimes u$.

The Euler-Poincar\'{e} principle implies 
\[ \Lambda^*(E) = H(E) \]
where $\Lambda^* (E)=\oplus_{n\geq 0} (-1)^n \Lambda^n (E)$ is the alternating sum of exterior powers of $E$ and $H(E)=\oplus_{n\geq 0} (-1)^n H_n(E)$ is the alternating sum of homology groups of $E$. Both sides of this identity are graded by the root lattice $II_{9,1}$ of $G$ and the homogeneous subspaces are finite dimensional. The homology groups $H_n(E)$ can be calculated in the same way as for Kac-Moody algebras (cf. \cite{B}). The result is that $H_n(E)$ is the subspace of $\Lambda^n (E)$ spanned by the homogeneous vectors of $\Lambda^n (E)$ of degree $\al \in II_{9,1}$ with $\al^2=0$. We can work out the homology groups of $G$ explicitly because we know the simple roots. Denote the subspace of $H(E)$ with degree $\al \in II_{9,1}$ by $H(E)_{\al}$ and analogous for $E$. Let $\lambda$ be a primitive norm zero vector in $II_{9,1}^+$. Then 
${\mathbb R} \oplus \oplus_{n>0} H(E)_{n \lambda} 
                   = \Lambda^*(\oplus_{n>0}E_{n\lambda})$. 
The denominator identity of the fake monster superalgebra now follows easily by calculating the character on both sides of $\Lambda^*(E) = H(E)$.

More generally if $g$ is an automorphism of $G$ then $g$ commutes with the derivations $d_i : \Lambda^i(E) \rightarrow \Lambda^{i-1}(E)$ and we get a sequence\[ \ldots \rightarrow g(\Lambda^2(E)) \rightarrow g(\Lambda^1(E)) 
          \rightarrow g(\Lambda^0(E))  \rightarrow 0 \, .  \] 
$g$ induces an isomorphism between the homology groups of the two complexes above. Hence we can apply $g$ to both sides of the equation $\Lambda^*(E) = H(E)$ and take the trace. This gives a twisted denominator identity. It depends only on the conjugacy class of $g$ in the automorphism group of $G$. 

Let $u$ be an element of $2.Aut(E_8)$ of odd order $N$.

The lattice $E_8$ has a unique central extension $\hat{E}_8$ by $ \{1,-1\}$ such that the commutator of any inverse images of $r,v$ in $E_8$ is $(-1)^{(r,v)}$. Since $\rho_V(u)$ has odd order it has a lift to $Aut (\hat{E}_8)=2^8.Aut (E_8)$ such that $\rho_V(u)^n$ fixes all elements of $\hat{E}_8$ which are in the inverse image of the vectors of $E_8$ fixed by $\rho_V(u)^n$ (cf. Lemma 12.1 in \cite{B}). We define $g$ as the automorphism of $G$ induced by this lift. 

Let $E_8^u$ be the sublattice of $E_8$ fixed by $\rho_V(u)$. The natural projection 
$\pi :{\mathbb R} \otimes E_8 \rightarrow 
      {\mathbb R} \otimes E_8^u$ 
maps $E_8$ onto the dual lattice $E_8 ^{u*}$ because $E_8$ is unimodular. We define the Lorentzian lattice $L=E_8^u\oplus II_{1,1}$ with dual lattice $L^* = E_8 ^{u*} \oplus II_{1,1}$ and denote the closures of the canonical positive cones by $L^+$ and $L^{*+}$. 
The fake monster superalgebra has a natural $L^*$-grading.   
For $\al=(r^*,a)\in L^*$ we define  
\[ \tilde{E}_{0,\al} = \oplus_{\pi(r)=r^*} E_{0,(r,a)}  \]
and analogous for $\tilde{E}_{1,\al}$. 

Let $\varepsilon_i$ and $\sigma_i$ denote the eigenvalues of $\rho_V(u)$ and $\rho_L(u)$. 

Then we have 
\begin{th1}
The twisted denominator identity corresponding to g is given by 
\[ \prod_{\al \in L^{*+}}
   \frac{ (1-e^{\al})^{ {mult}_0({\al}) }}
        { (1+e^{\al})^{ {mult}_1({\al}) }} =
    1 + \sum a(\lambda)e^{\lambda}  \]
where
\beann 
\mbox{mult}_0({\al}) &=& 
\sum_{ds|(({\al},L),N)} \frac{\mu(s)}{ds}\, tr(g^d|\tilde{E}_{0,{\al}/ds}) \\
\mbox{mult}_1({\al}) &=& 
\sum_{ds|(({\al},L),N)} \frac{\mu(s)}{ds}\, tr(g^d|\tilde{E}_{1,{\al}/ds})
\eeann
and $a(\lambda)$ is the coefficient of $q^m$ in
\[ \prod_{n\geq 1} \prod_{\: 1\leq i\leq 8}
   \frac{(1-\varepsilon_i q^n)}
      {(1+\sigma_i q^n)}      \]
if $\lambda$ is $m$ times a primitive norm zero vector in $L^+$ and zero else. 
\end{th1}
{\em Proof:} We consider both sides of $\Lambda^*(E) = H(E)$ as $L^*$-graded $g$-modules. $\Lambda^*(E)$ is isomorphic to $\Lambda^*(E_0) \otimes S^*(E_1)$ if we forget the superstructure on $E_1$. First we calculate the trace of $g$ on $S^*(E_1)$. For that we recall some formulas. Let $V$ be a finite dimensional vector space. Then we have
\[  \sum _{n\geq 0}(-1)^n (dim \, S^n(V)) q^n = (1+q)^{-dim V} 
    = exp \, \Big\{  \sum_{n>0}(-1)^n (dim \,V) q^n/n  \Big\}  \, . \]
The second equality can be proven by taking logarithms and using the expansion $log (1+q) =- \sum_{n>0}  (-1)^n q^n/n$.
This can be generalized to 
\[    \sum _{n\geq 0}(-1)^n \,tr(g|S^n(V)) q^n 
     = exp \, \Big\{ \sum_{n>0} (-1)^n tr(g^n|V) q^n/n  \Big\}  \, . \]
Another formula we need is
$S^*(V_1 \oplus V_2) = S^*(V_1) \otimes S^*(V_2)$. Using these two formulas we find that the trace of $g$ on $S^*(E_1)$ is given by 
\[  exp \, \Big\{ 
      \sum_{\bt \in L^{*+}} \sum_{n>0} (-1)^n \,tr(g^n|\tilde{E}_{1,\bt} ) 
              e^{n\bt}/n                        \Big\} \, . \]
We want to express the trace in the form 
\[ \prod_{\bt \in L^{*+}}(1+e^{\bt})^{ - mult_1({\bt}) }=
   \mbox{\it exp} \,
   \Big\{ \sum_{\bt \in L^{*+}} \sum_{n>0} (-1)^n mult_1(\bt)e^{n\bt}/n 
       \Big\}  \, . \]
Taking logarithms and comparing coefficients at $e^{\al }$ this implies
\[  \sum_{\bt \in L^{*+} \atop n\bt   = \al } 
           (-1)^n \,tr(g^n|\tilde{E}_{1,\bt} )/n =
    \sum_{\bt \in L^{*+} \atop n\bt  = \al } 
           (-1)^n mult_1(\bt) /n   \, . \]
The vector $\al \in L^*$ is a positive multiple $m$ of a primitive vector in $L^*$. This vector generates a lattice that we can identify with ${\mathbb Z}$. Then the right hand side of the last equation is the convolution product of the arithmetic functions $h(n)=(-1)^n/n$ and $mult_1(n)$. Hence 
\[ 
mult_1(\al)= \sum_{ds|m} h^{* -1}(s) (-1)^d\, tr(g^d|\tilde{E}_{1,\al /ds})/d  
\, . \]
Let $\mu$ be the M\"obiusfunction and $f(n)$ the arithmetic function which is zero if $n$ contains an odd square and $(-1)^k / p_1\ldots p_k$, where the $p_i$ are the different primes dividing $n$, else. Then convolution inverse of $h(n)$ is given by 
\[ h^{* -1}(n) = \left\{ \begin{array}{cl}
                           (\mu h) (n)  & \quad n \quad \mbox{odd} \\
                            f(n)        & \quad n \quad \mbox{even}
                         \end{array}
                 \right. \] 
$g$ has the same order as $u$ because $N$ is odd. The trace $tr(g^d|\tilde{E}_{1,\al})$ depends only on $\al$ and $(d,N)$. This implies that the sum in the expression for $mult_1(\al)$ extends only over $ds|(m,N)$. It is easy to see that $m$ is equal to the highest common factor $(\al,L)$ of the numbers 
$(\al,\bt),\, \bt \in L$. Hence we get the following formula 
\[ 
mult_1(\al) = \sum_{ds|((\al,L),N)} \mu(s)\, tr(g^d|\tilde{E}_{1,\al/ds})/ds
\, . \]
Note that this formula is wrong for even $N$. This is another reason why we restrict to elements in $2.Aut(E_8)$ of odd order.

If we express the contributions of $\Lambda^*(E_0)$ to the trace in the form 
\[ \prod_{\al \in L^{*+}}(1-e^{\al})^{ mult_0(\al) } \]
then we find 
\[ mult_0 (\al)= \sum_{ds|((\al,L),N)} h^{* -1}(s)\, 
                                     tr(g^d|\tilde{E}_{0,\al/ds})/d \]
with $h(n)=1/n$. This function is strongly multiplicative so that its convolution inverse is given by $\mu h$. An argument as above then gives the expression for $mult_0(\al)$.  

Next we calculate the trace of $g$ on $H(E)$. We have 
\[ H(E) = \oplus_{\al \in L^{*+}} \tilde{H}(E)_{\al} \] with 
\[ \tilde{H}(E)_{\al}=\oplus_{\pi(r)=r^*} H(E)_{(r,a)} \]
for $\al=(r^*,a)$ in $L^{*+}$.
Clearly $tr(g|\tilde{H}(E)_{\al})$ is zero if $\al$ is not in $L$ and 
$tr(g|\tilde{H}(E)_{\al})=tr(g|H(E)_{\al})$ for $\al\in L$.
We can also restrict to $\al^2=0$. Let $\al\in L^+ \subset II_{9,1}^+$ be m times a primitive norm zero vector $\lambda$ in $II_{9,1}^+$. Then $\lambda$ is also a primitive vector in $L$. $H(E)_{\al}$ is determined by
\[  {\mathbb R} \oplus \oplus_{n>0} H(E)_{n \lambda} 
    = \Lambda^*(\oplus_{n>0}E_{0,n\lambda}) \otimes 
      S^*(\oplus_{n>0}E_{1,n\lambda}) \, .               \]
g acts by $\rho_V(u)$ on $E_{0,n\lambda}\cong {\mathbb R}\otimes E_8$ and by 
$\rho_L(u)$ on $E_{1,n\lambda}\cong C$ (cf. section 2). 
Going over to complexifications this implies that $tr(g|H(E)_{\al})$ is given by the coefficient of $q^m$ in 
\[ 1 + \sum_{n \geq 1} tr(g|H(E)_{n \lambda}) q^n = 
                  \prod_{n\geq 1} \prod_{\: 1\leq i\leq 8}
                  \frac{(1-\varepsilon_i q^n)}{(1+\sigma_i q^n)} \, . \]
This finishes the proof of the theorem.

For $u=1$ the theorem gives the denominator identity of the fake monster superalgebra. 

If the multiplicities are nonnegative integers then the twisted denominator identity is the untwisted denominator identity of a generalized Kac-Moody superalgebra with the following properties. The root lattice is $L^*$ or a sublattice thereof. The algebra has no real roots so that the Weyl group is trivial. The multiplicities of a root $\al$ are given by ${mult}_0({\al})$ and ${mult}_1({\al})$. There are no real simple roots and the imaginary  simple roots are the norm zero vectors in $L^{*+}$. The Weyl vector is zero. We describe the multiplicities of the simple roots. Let $\al = n \lambda$ be a simple root where $\lambda$ is a primitive vector in $L^+$. Suppose that $\rho_V(u)$ and $\rho_L(u)$ have cycle shapes $a_1^{b_1}\ldots a_k^{b_k}$ and $c_1^{d_1}\ldots c_l^{d_l}$. Then the multiplicities of $\al$ as even and odd root are given by 
$\sum_{a_k | n} b_k$ and $\sum_{c_k | n} d_k$. 

The formulas for the multiplicities simplify if $u$ has in addition prime order.
In this case we only need to calculate the trace of $g$ on $\tilde{E}_{0,\al}$ and $\tilde{E}_{1,\al}$ with $\al \in L^*$ and the dimensions of these spaces. We find
\begin{pth1} \label{hel}
For $\al \in L^*$ the trace $tr(g|\tilde{E}_{0,\al})$ is given by the coefficient of $q^{(1-\al^2)/2}$ in 
\[ \frac{1}{2}\left( 
    \prod_{n \geq 1} \prod_{\: 1 \leq i \leq 8} 
      \frac{(1+\varepsilon_i q^{n-1/2})}
           {(1-\varepsilon_iq^n)} 
  - \prod_{n \geq 1} \prod_{\: 1 \leq i \leq 8}
      \frac{(1-\varepsilon_i q^{n-1/2})}
           {(1-\varepsilon_iq^n)}
              \right) \] 
if $\al \in L$ and zero else. 

Let $tr (\rho_L(u))=tr (\rho_R(u))$. Then $tr(g|\tilde{E}_{1,\al})$ is the coefficient 
of $q^{(1-\al^2)/2}$ in
\[ tr (\rho_L(u)) \: q^{1/2}
 \prod_{n \geq 1} \prod_{\: 1 \leq i \leq 8}
    \frac{(1+\varepsilon_iq^n)}
         {(1-\varepsilon_iq^n)}
\]
if $\al \in L$ and zero else. 
\end{pth1}
{\em Proof:} Clearly the trace of $g$ is zero on $\tilde{E}_{0,\al}$ if $\al$ is not in $L$. For $\al \in L$ we have $tr(g|\tilde{E}_{0,\al})=tr(g|E_{0,\al})$. Write $\al = (r,a)$. Then $E_{0,\al}$ is isomorphic as $g$-module to the subspace of $V_{0,(1-a^2)/2}$ of degree $r$ (cf. section 2). This space is generated by products of fermionic and bosonic oscillators and $e^r$. The sum of the $L_0$-contribution of the oscillators and the $L_0$-contribution of $e^r$ is $\frac{1}{2}-\frac{1}{2}a^2$. The vector $e^r$ has $L_0$-eigenvalue $\frac{1}{2}r^2$ so that the $L_0$-contribution of the oscillators is $ \frac{1}{2}-\frac{1}{2}a^2-\frac{1}{2} r^2 = \frac{1}{2}-\frac{1}{2}\al^2$. Now we go over to complexifications and choose a basis of ${\mathbb C}\otimes E_8$ in which $\rho_V(u)$ is diagonal. Then the trace of $g$ on $E_{0,\al}$ is given by the coefficient of $q^{(1-\al^2)/2}$ in 
\[ \prod_{n \geq 1} \prod_{\: 1 \leq i \leq 8} 
      \frac{(1+\varepsilon_i q^{n-1/2})}
           {(1-\varepsilon_iq^n)} \,. \] 
Since we only need the half integral exponents of $q$ in this expression we can subtract the integral exponents. This proves the first statement. 

The argument for the second statement is similar. Note that there are two types of ground states on which $u$ may act differently. We avoid this problem by assuming that $u$ has in both representations the same trace.
This proves the proposition.
\begin{pth1} \label{de}
Let $\al=(r^*,a)\in L^*$ and let $r^{\bot *}\in E_8^{u \bot *}$ such that $r^* + r^{\bot *} \in E_8$. Then the dimension of $\tilde{E}_{0,\al}$ and $\tilde{E}_{1,\al}$ is given by the coefficient of $q^{(1-\al^2)/2}$ in 
\[ 8 q^{1/2} \frac{\eta(q^2)^8}{\eta(q)^{16}}\:
   \theta_{r^{\bot *} +E_8^{u \bot }}(q) \]
where $\theta_{r^{\bot *}+E_8^{u \bot }}(q)$ is the theta function of the translated lattice $r^{\bot *}+E_8^{u \bot }$. 
\end{pth1}
{\em Proof:} Clearly $\tilde{E}_{0,\al}$ and $\tilde{E}_{1,\al}$ have the same dimension. The inverse image of $r^*$ in $E_8$ under $\pi$ is $r^*+(r^{\bot *}+E_8^{u \bot })$. Hence $\tilde{E}_{1,\al}=\oplus_{\pi(r)=r^*} E_{1,(r,a)}$ is isomorphic to the direct sum of the $V_{1,(1-a^2)/2}(r^*+s)$ where $s$ is in the translated lattice $r^{\bot *}+E_8^{u \bot}$. As above the sum of the $L_0$-contribution of the bosonic and fermionic oscillators and $\frac{1}{2}s^2$ is $\frac{1}{2}-\frac{1}{2}\al^2$. The proposition now follows from simple counting.

We will use the shorter notation $\theta_{r^{\bot *}}(q)$ in the following.

\section{Two supersymmetric algebras}

In this section we calculate explicitly twisted denominator identities corresponding to elements in $2.Aut(E_8)$ of order 3 and 7. These identities are the untwisted denominator identities of 2 supersymmetric generalized Kac-Moody superalgebras.
 
We choose an orthonormal basis $\{e_0,e_1,\ldots,e_7 \}$ of ${\mathbb R}^8$ and embed the lattice $E_8$ as the set of points $\sum m_i e_i$ where all $m_i$ are in $\mathbb Z$ or all $m_i$ are in $\mathbb Z+\frac{1}{2}$ and $\sum m_i$ is even. In these coordinates the automorphism group of $E_8$ is generated by the permutations of the coordinates, even sign changes and an involution generated by the Hadamard matrix. We also identify ${\mathbb R}^8$ with the alternative algebra ${\mathbb O}$ of octonions by defining $e_0$ as the identity and $e_i e_j= a_{ijk}e_k-\delta_{ij}1$ for $1\leq i,j \leq 7$ where $a_{ijk}$ is the totally antisymmetric tensor with 
$a_{ijk}=1$ for $ijk=123, 154, 264, 374, 176, 257, 365$.

The element $u=\frac{1}{4}1(e_2-e_3)1(e_1-e_2)1(e_6-e_7)1(e_5-e_6)$ in $2.Aut(E_8)$ has order $3$. It is easy to check that the transformations corresponding to $\rho_V(u),\rho_L(u)$ and $\rho_R(u)$ are all equal and 
\[ \renewcommand{\arraystretch}{1.3} 
\begin{array}{llll}  
\rho_V(u)1=1 & \rho_V(u)e_1 = e_3 & \rho_V(u)e_2 = e_1 & \rho_V(u)e_3 = e_2 \\ 
\rho_V(u)e_4=e_4 & \rho_V(u)e_5 = e_7 & \rho_V(u)e_6 = e_5 
                                                   & \rho_V(u)e_7 = e_6 \, . 
\end{array} \]
Hence $\rho_V(u),\rho_L(u)$ and $\rho_R(u)$ all have cycle shape $1^23^2$. 

The following proposition collects some results on $E^u_8$.
\begin{p1}
The sublattice $E^u_8$ of $E_8$ fixed by $\rho_V(u)$ is the 4 dimensional lattice with elements $(m_1,m_2,m_3,m_4)$ where all $m_i$ are in $\mathbb Z$ or all $m_i$ are in $\mathbb Z+\frac{1}{2}$ and $\sum m_i$ is even. The norm is $(m_1,m_2,m_3,m_4)^2=m_1^2+m_2^2+3m_3^2+3m_4^2$. $E^u_8$ has determinant $3^2$ and the quotient $E_8^{u*}/E_8^u$ is ${\mathbb Z}_3\times {\mathbb Z}_3$. The level of $E_8^u$ is $3$ so that $3E_8^{u*}$ is a sublattice of $E^u_8$. The orthogonal complement of $E^u_8$ in $E_8$ is $A_2 \oplus A_2$. 
\end{p1}

We recall some results about modular forms. The congruence subgroup of $SL_2({\mathbb Z})$ of level $3$ is defined as  
$\Gamma(3)=\{ \big( {a \atop c} {b \atop d} \big) \in SL_2({\mathbb Z}) \, |\,
   \big( {a \atop c} {b \atop d} \big) = \big( {1 \atop 0} {0 \atop 1} \big)
   \mbox{ mod } 3 \, \}$.
The vector space of modular forms for $\Gamma(3)$ with even positive weight $n$ has dimension $n+1$. A modular form in this vector space is zero if and only if the coefficients of $q^0,q^{1/3},\ldots,q^{n/3}$ in its Fourier expansion are zero. 

For $r\in A_2^*\oplus A_2^*$ we define $\delta(r)=1$ if $r\in A_2 \oplus A_2$ and $0$ else. We have
\begin{pp1} \label{tta}
Let $r\in A_2^*\oplus A_2^*$. Then
\[ \theta_{r+ A_2\oplus A_2}(q) =
\frac{1}{4} \, \frac{\eta(q)^{12}}{\eta(q^2)^6} 
\left\{
  \frac{\eta(q^6)^2}{\eta(q^3)^4}\delta(r) 
+ \sum_{j=0}^2 \varepsilon^{-3jr^2/2} 
             \frac{ \eta\!\left( (\varepsilon^jq^{1/3})^2 \right)^2 }
                  { \eta\!\left( \varepsilon^jq^{1/3} \right)^4    } 
\right\} \]
where $\varepsilon = e^{2\pi i/3}$.
\end{pp1}
{\em Proof:} $A_2\oplus A_2$ has dimension $4$ and level $3$ so that $\theta_{r+ A_2\oplus A_2}$ is a modular form for $\Gamma(3)$ of weight $2$. The same holds for the right hand side of the formula which can be shown by calculating the transformations under a set of generators. The proposition follows from comparing the coefficients at $q^0,q^{1/3}$ and $q^{2/3}$.   

The next identity is a twisted version of Jacobi's identity.
\begin{pp1} \label{hnel}
Let $|q|<1$. Then 
\[ \frac{1}{2q^{1/2}}
   \left\{
    \prod_{n \geq 1}(1+q^{3n-3/2})^2(1+q^{n-1/2})^2
   -\prod_{n \geq 1}(1-q^{3n-3/2})^2(1-q^{n-1/2})^2
   \right\} \]
\[   = 2 \prod_{n \geq 1} (1+q^{3n})^2(1+q^n)^2 \, . \]
\end{pp1}
{\em Proof:} This is an identity between modular forms and therefore can be proven by comparing sufficiently many coefficients in their Fourier expansions.  
 
The proposition implies that the generating functions for $tr(g|\tilde{E}_{0,\al})$ and \linebreak $tr(g|\tilde{E}_{1,\al})$ are equal.

Define $c(n)$ by 
\[  \sum_{n\geq 0} c(n)q^n 
   = 2\, \prod_{n \geq 1} \frac{(1+q^{3n})^2(1+q^n)^2 }{(1-q^{3n})^2(1-q^n)^2 }
   = 2\, \frac{\eta(q^6)^2\eta(q^2)^2}{\eta(q^3)^4\eta(q)^4} \]
\[   = 2 + 8q + 24q^2 + 72q^3 + 184q^4 + 432q^5 + 984q^6 + 2112q^7 + \ldots \]
Let $g$ be the automorphism induced by $u$. Then we have
\begin{thp1}
The twisted denominator identity corresponding to $g$ is
\[ \prod_{\al\in L^{+}}
     \frac{ (1-e^{\al})^{ c(-\al^2/2) }}
          { (1+e^{\al})^{ c(-\al^2/2) }}
   \prod_{\al \in L^{+}\cap 3L^*}
     \frac{ (1-e^{\al})^{ c(-\al^2/6) }}
          { (1+e^{\al})^{ c(-\al^2/6) }}
=  1 + \sum a(\lambda)e^{\lambda}  \]
where $a(\lambda)$ is the coefficient of $q^n$ in 
\[  \prod_{n \geq 1}\frac{(1-q^{3n})^2(1-q^n)^2}{(1+q^{3n})^2(1+q^n)^2}
    = 1 - 4q + 4q^2 - 4q^3 + 20q^4 - 24q^5 + 4q^6 - \ldots  \]
if $\lambda$ is $n$ times a primitive norm zero vector in $L^+$ and zero else. 
\end{thp1}
{\em Proof:} Recall that $3L^*\subset L$ because $L$ has level 3. Furthermore 
$3$ divides $(\al,L)$ if and only if $\al \in 3L^*$. 
We consider now 4 cases. \\
$\al\notin L$. 
Then $\al\notin 3L^*$ and by Proposition \ref{hel} the multiplicities $mult_0(\al)$ and $mult_1(\al)$ are zero.\\
$\al \in L$ and $\al\notin 3L^*$. 
Then $mult_0(\al)=tr(g|\tilde{E}_{0,\al})$ and $mult_1(\al)=tr(g|\tilde{E}_{1,\al})$. Using Propositions \ref{hel} and \ref{hnel} we find $mult_0(\al)=mult_1(\al)=c(-\al^2/2)$.\\
$\al\in 3L^*$ and $\al \notin 3L$. 
Then $mult_0(\al)=tr(g|\tilde{E}_{0,\al})+tr(1|\tilde{E}_{0,\al/3})/3$. The first term gives $c(-\al^2/2)$. Write $\al/3=(r^*,a)$ where $r^*$ is in $E_8^{u*}$ but not in $E_8^{u}$ and choose $r^{*\bot}$ as in Proposition \ref{de}. Then the dimension of $\tilde{E}_{0,\al/3}$ 
is the coefficient of $q^{-\al^2/18}$ in 
$8 \eta(q^2)^8 \theta_{r^{\bot *}}(q)/{\eta(q)^{16}}$ 
or equivalently the coefficient of $q^{-\al^2/6}$ in 
$8 \eta(q^6)^8 \theta_{r^{\bot *}}(q^3)/ \eta(q^3)^{16}$. 
Note that $\al^2\in 6 {\mathbb Z}$. We have 
$\al^2/9 = r^{*2} = - r^{\bot*2}$ mod $2$ and 
$\al^2/6 = - 3 r^{\bot*2}/2$ mod $3$ so that 
by Proposition \ref{tta}
\[ \theta_{r^{\bot *}}(q^3) = 
      \frac{1}{4} \, \frac{\eta(q^3)^{12}}{\eta(q^6)^6} 
      \sum_{j=0}^2 \varepsilon^{j\al^2/6} 
             \, \frac{ \eta\left( (\varepsilon^j q)^2 \right)^2 }
                  { \eta\left( \varepsilon^j q \right)^4    } 
\, . \]
The coefficient of $q^{-\al^2/6}$ in $8 \eta(q^6)^8 \theta_{r^{\bot *}}(q^3)/ \eta(q^3)^{16}$ is equal to the coefficient of $q^{-\al^2/6}$ in 
\[ 6 \, \frac{\eta(q^6)^2\eta(q^2)^2}{\eta(q^3)^4\eta(q)^4}\, . \]
This shows that $mult_0(\al) = c(-\al^2/2) + c(-\al^2/6)$. 
The result for $mult_1(\al)$ is clear. \\
$\al \in 3L$. Then 
$mult_0(\al)=tr(g|\tilde{E}_{0,\al})-tr(g|\tilde{E}_{0,\al/3})/3 
                          + tr(1|\tilde{E}_{0,\al/3})/3$. 
Here we have an additional term in $\theta_{r^{\bot *}}$ which cancels exactly with the term from $tr(g|\tilde{E}_{0,\al/3})$ so that we get the same result for $mult_0(\al)$ as in the case before. The result for $mult_1(r)$ is again clear. \\
This proves the theorem.
 
Since the multiplicities are all nonnegative integers there is a generalized Kac-Moody superalgebra whose denominator identity is the identity given in the theorem.
\begin{cp1}
There is a generalized Kac-Moody superalgebra with root lattice $L$ and root multiplicities given by
\[ mult_0(\al) = mult_1(\al) = c(-\al^2/2) 
                               \qquad \al\in L,\: \al \notin 3L^* \]
and 
\[ mult_0(\al) = mult_1(\al) = c(-\al^2/2)+ c(-\al^2/6)  
                                           \qquad \al\in 3L^*\,. \]
The simple roots of are the norm zero vectors in $L^+$. Their multiplicities as even and odd roots are equal. Let $\lambda$ be a simple root. Then $mult_0(\lambda)=4$ if $\lambda$ is 3n times a primitive vector in $L^+$ and $mult_0(\lambda)=2$ else.
\end{cp1}
As the fake monster superalgebra this algebra is supersymmetric and has no real roots. The Weyl group is trivial and the Weyl vector is zero. The supersymmetry is a consequence of the twisted Jacobi identity in Proposition \ref{hnel}. In a forthcoming paper we show that the denominator function of this algebra defines an automorphic form for a discrete subgroup of $O_{6,2}({\mathbb R})$ of weight $2$.  

Now we consider the next example. The analysis is similar to the above one. The element $u=\frac{1}{8}1(e_6-e_7)1(e_5-e_6)1(e_4-e_5)1(e_3-e_4)1(e_2-e_3)1(e_1-e_2)$ in $2.Aut(E_8)$ has order $7$. The transformations corresponding to $\rho_V(u),\rho_L(u)$ and $\rho_R(u)$ are all equal and 
\[ \renewcommand{\arraystretch}{1.3} 
\begin{array}{llll}  
\rho_V(u)1=1 & \rho_V(u)e_1 = e_7 & \rho_V(u)e_2 = e_1 & \rho_V(u)e_3 = e_2 \\ 
\rho_V(u)e_4=e_3 & \rho_V(u)e_5 = e_4 & \rho_V(u)e_6 = e_5 
                                                   & \rho_V(u)e_7 = e_6 \, . 
\end{array} \]
Hence $\rho_V(u),\rho_L(u)$ and $\rho_R(u)$ have cycle shape $1^17^1$. We have
\begin{pp1}
The sublattice $E^u_8$ of $E_8$ fixed by $\rho_V(u)$ is the 2 dimensional lattice with elements $(m_1,m_2)$, where either $m_1$ and $m_2$ are in $\mathbb Z$ and $m_1+m_2$ is even or $m_1$ and $m_2$ are in $\mathbb Z+\frac{1}{2}$ and $m_1+m_2$ is odd, and norm $(m_1,m_2)^2=m_1^2+7m_2^2$. The quotient $E_8^{u*}/E_8^u$ is ${\mathbb Z}_7$ and $E^u_8$ has level $7$. The orthogonal complement of $E^u_8$ in $E_8$ is isomorphic to $A_6$. 
\end{pp1}
The theta function $\theta_{r+ A_6}$  of a coset $r+A_6$ of $A_6$ in its dual depends only on $r^2$ mod $2$.
\begin{pp1} 
Let $r\in A_6^*$. Then
\[ \theta_{r+ A_6}(q) =
\frac{1}{8}\, \frac{\eta(q)^{14}}{\eta(q^2)^7} 
\left\{
  \frac{\eta(q^{14})}{\eta(q^7)^2}\delta(r) 
+ \sum_{j=0}^6 \varepsilon^{-7jr^2/2} 
             \frac{ \eta\left( (\varepsilon^jq^{1/7})^2 \right)^2 }
                  { \eta\left( \varepsilon^jq^{1/7} \right)^4    } 
\right\} \]
where $\varepsilon = e^{2\pi i/7}$.
\end{pp1}
The following supersymmetry relation holds.
\begin{pp1}
Let $|q|<1$. Then 
\[  \frac{1}{2q^{1/2}}
    \left\{ 
      \prod_{n \geq 1}(1+q^{7n-7/2})(1+q^{n-1/2})
     -\prod_{n \geq 1}(1-q^{7n-7/2})(1-q^{n-1/2})
    \right\} \]
\[  = \prod_{n \geq 1} (1+q^{7n})(1+q^n)  \]
\end{pp1}
Here we define the numbers $c(n)$ by 
\[  \sum_{n\geq 0} c(n)q^n 
   =  \prod_{n \geq 1} \frac{(1+q^{7n})(1+q^n)}{(1-q^{7n})(1-q^n)}
   =  \frac{\eta(q^{14})\eta(q^2)}{\eta(q^7)^2\eta(q)^2} \]
\[   =  1 + 2q + 4q^2 + 8q^3 + 14q^4 + 24q^5 + 40q^6 + 66q^7 + \ldots \]
Write $g$ for the automorphism induced by $u$. Then 
\begin{thp1}
The twisted denominator identity corresponding to $u$ is
\[ \prod_{\al\in L^{+}}
   \frac{ (1-e^{\al})^{ c(-\al^2/2) } }
        { (1+e^{\al})^{ c(-\al^2/2) } } 
   \prod_{\al \in L^{+}\cap 7L^*}
   \frac{ (1-e^{\al})^{ c(-\al^2/14) } }
        { (1+e^{\al})^{ c(-\al^2/14) } }
=  1 + \sum a(\lambda)e^{\lambda}  \]
where $a(\lambda)$ is the coefficient of $q^n$ in 
\[  \prod_{n \geq 1}\frac{(1-q^{7n})(1-q^n)}{(1+q^{7n})(1+q^n)}
    =  1 - 2q + 2q^4 - 2q^7 + 4q^8 - 2q^9 - 4q^{11} + 6q^{16} - \ldots  \]
if $\lambda$ is $n$ times a primitive norm zero vector in $L^+$ and zero else. 
\end{thp1}
Again the multiplicities are all nonnegative integers and we have 
\begin{cp1}
There is a generalized Kac-Moody superalgebra with root lattice $L$ and root multiplicities given by
\[ mult_0(\al) = mult_1(\al) = c(-\al^2/2) 
                                   \qquad \al \in L,\: \al \notin 7L^* \]
and 
\[ mult_0(\al) = mult_1(\al) = c(-\al^2/2)+ c(-\al^2/14)  
                                           \qquad \al \in 7L^*\,. \]
The simple roots of are the norm zero vectors in $L^+$. Let $\lambda$ be a simple root. Then $mult_0(\lambda)=mult_1(\lambda)=2$ if $\lambda$ is 7n times a primitive vector in $L^+$ and $mult_0(\lambda)=mult_1(\lambda)=1$ else. 
\end{cp1}
The denominator identity of this algebra is given by the identity in the above theorem. It defines an automorphic form for a subgroup of $O_{4,2}({\mathbb R})$ of weight $1$.

\section*{Acknowledgments}

I thank R. E. Borcherds for stimulating discussions.

\end{document}